\documentclass[11pt]{article}
\usepackage{graphicx}
\def\rsq{\hfill\rule{1mm}{3mm}}

\newtheorem{guess}{\bf Theorem}

\newtheorem{lemma}{\bf Lemma}

\begin{document}
\begin{center}
{\Large A note on 3-colorable plane graphs without 5- and 7-cycles
\footnote{\footnotesize Supported partially by NSFC 10371055}}

\vskip 15pt

{Baogang Xu\footnote{\footnotesize email: baogxu@njnu.edu.cn}}

{\small School of Mathematics and Computer Science, Nanjing Normal
University, 122 Ninghai Road}

{\small Nanjing, 210097, PR China}

\end{center}

\vspace{0.5cm} 
\begin{abstract}
In \cite{bgrs0}, Borodin {\em et al} figured out a gap of
\cite{bxu1}, and gave a new proof with the similar technique. The
purpose of this note is to fix the gap of \cite{bxu1} by slightly
revising the definition of {\em special faces}, and adding a few
lines of explanation in the proofs (new added text are all in
black font).
\begin{flushleft}
{\em Key words and phrases: plane graph, cycle, coloring

}
\end{flushleft}
\end{abstract}

In \cite{bgrs}, Borodin {\it et al}  proved that every plane graph
$G$ without cycles of length from 4 to 7 is 3-colorable that
provides a new upper bound to Steinberg's conjecture (see
\cite{a11} p.229). In \cite{ovbar}, Borodin and Raspaud proved
that every plane graph with neither 5-cycles nor triangles of
distance less than four is 3-colorable, and they conjectured that
every plane graph with neither 5-cycles nor adjacent triangles is
3-colorable, where the distance between triangles is the length of
the shortest path between vertices of different triangles, and two
triangles are said to be adjacent if they have an edge in common.
In \cite{bxu2}, Xu improved Borodin and Raspaud's result by
showing that every plane graph with neither 5-cycles nor triangles
of distance less than three is 3-colorable.

In this note, it is proved that every plane graph without 5- and
7-cycles and without adjacent triangles is 3-colorable. This
improves the result of \cite{bgrs}, and offers a partial solution
for Borodin and Raspaud's conjecture \cite{ovbar}.

Let $G=(V, E, F)$ be a plane graph, where $V, E$ and $F$ denote
the sets of vertices, edges and faces of $G$ respectively. The
neighbor set and degree of a vertex $v$ are denoted by $N(v)$ and
$d(v)$, respectively. Let $f$ be a face of $G$. We use $b(f),
V(f)$ and $N(f)$ to denote the boundary of $f$, the set of
vertices on $b(f)$, and the set of faces adjacent to $f$
respectively. The degree of $f$, denoted by $d(f)$, is the length
of the facial walk of $f$. A $k$-vertex ($k$-face) is a vertex
(face) of degree $k$.

Let $C$ be a cycle of $G$. We use $int(C)$ and $ext(C)$ to denote
the sets of vertices located inside and outside $C$, respectively.
$C$ is called a {\it separating} cycle if both $int(C)\neq
\emptyset$ and $ext(C)\neq \emptyset$, and is called a {\it facial
cycle} otherwise. For convenience, we still use $C$ to denote the
set of vertices of $C$.

{\bf Let $f$ be an 11-face bounded by a cycle  $C=u_1u_2u_3\ldots
u_{11}u_1$. A 4-cycle $u_1u_2u_3vu_1$ is called an {\em ear} of
$f$ if $v\not\in C$. The graph $G_1$, obtained from $G$ by
removing $u_2$ and all the vertices in $int(u_1u_2u_3vu_1)$,  is
called an {\em ear-reduction} of $G$ on $f$. Since $u_1vu_3\ldots
u_{11}u_1$ is still an 11-cycle bounding a face, say $f_1$, in
$G_1$, if $f_1$ has an ear, we may make an ear-reduction to $G_1$
on $f_1$ and get a new graph $G_2$ and an 11-face $f_2$ bounded by
a cycle in $G_2$. Continue this procedure, we get a sequence of
graphs $G, G_1, G_2, \ldots$, and a sequence of 11-faces $f, f_1,
f_2, \ldots $, such that $f_i$ is an 11-face in $G_i$. Each of
these 11-faces is called a {\em collapse} of $f$.

An 11-face $f$ of $G$ is called a {\it special face} if the
following hold: (1) $b(f)$ is a cycle; (2) $f$ is adjacent to a
triangle sharing only one edge with $f$; and furthermore, for each
collapse $f'$ of $f$ and its corresponding graph $G'$: (3) every
vertex in $V(G')\setminus V(f')$ has at most two neighbors on
$b(f')$; and (4) for every edge $uv$ of $G'\setminus V(f')$,
$|N_{G'}(u)\cap V(f')|+|N_{G'}(v)\cap V(f')|\leq 3$.}

A vertex in $G\setminus V(f)$ that violates (3) is called a {\it
claw-center} of $b(f)$, and a pair of adjacent vertices in
$G\setminus V(f)$ that violates (4) is called a {\it
d-claw-center} of $b(f)$.

A separating 11-cycle $C$ is called a {\it special cycle} if in
$G\setminus ext(C)$, $C$ is the boundary of a special face. We use
${\cal G}$ to denote the set of plane graphs without 5- and
7-cycles and without adjacent triangles. Following is our main
theorem.

\renewcommand{\baselinestretch}{1}
\begin{guess}\label{th-1} Let $G$ be a graph in ${\cal G}$ that contains cycles of
length $4$ or $6$, $f$ an arbitrary face that is a special face,
or a $3$-face, or a $9$-face with $b(f)$ being a cycle. Then, any
$3$-coloring of $f$ can be extended to $G$.
\end{guess}
\renewcommand{\baselinestretch}{2}

As a corollary of Theorem \ref{th-1}, every plane graph in ${\cal
G}$ is 3-colorable. To see this, let $G$ be a plane graph in
${\cal G}$. By Gr\"{o}tzsch's theorem, we may assume that $G$
contains triangles. {\bf If $G$ contains neither 4-cycles nor
6-cycles, then by Theorem 1.2 of \cite{bgrs}, $G$ is 3-colorable}.
Otherwise, for an arbitrary triangle $T$, any 3-coloring of $T$
can be extended to $int(T)$ and $ext(T)$, that yields a 3-coloring
of $G$.

\vskip 10pt


\noindent {\bf Proof of Theorem \ref{th-1}.} Assume that $G$ is a
counterexample to Theorem \ref{th-1} with minimum
$\sigma(G)=|V(G)|+|E(G)|$. Without loss of generality, assume that
the unbounded face $f_o$ is a special face, or a 3-face or a
9-face with $b(f)$ being a cycle, such that a 3-coloring $\phi$ of
$f_o$ cannot be extended to $G$. Let $C=b(f_o)$ and let $p=|C|$.
Then, every vertex not in $C$ has degree at least 3.

By our choice of $G$, $f_o$ has no ears if $p=11$, and neither
4-cycle nor 6-cycle is adjacent to triangles. Since $G\setminus
int(C')$ is still in ${\cal G}$ for any separating cycle $C'$ of
$G$, {\bf either by the minimality of $G$ or by Theorem 1.2 of
\cite{bgrs} (this will be used frequently but implicitly)},

\renewcommand{\baselinestretch}{1}
\begin{lemma}\label{clm0}
$G$ contains neither special cycles, nor separating
$k$-cycles,$k=3,9$.
\end{lemma}

\begin{lemma}\label{2-connect}
$G$ is $2$-connected. That is, the boundary of every face of $G$
is a cycle.
\end{lemma}\renewcommand{\baselinestretch}{2}

\indent Interested readers may find the proof of Lemma
\ref{2-connect} in \cite{bgrs} (see that of Lemma 2.2).

Let $C'$ be a cycle of $G$, and $u$ and $v$ two vertices on $C'$.
We use $C'[u,v]$ to denote the path of $C'$ clockwisely from $u$
to $v$, and let $C'(u,v)=C'[u,v]\setminus \{u, v\}$. Unless
specified particularly, we always write a cycle on its vertices
sequence clockwisely.

\medskip

\renewcommand{\baselinestretch}{1}
\begin{lemma}\label{chordless}
$C$ is chordless.\end{lemma}
\renewcommand{\baselinestretch}{2}
{\bf Proof.} Assume to the contrary that $C$ has a chord $uv$. Let
$S_1=V(C(u,v))$, $S_2=V(C(v,u))$, and assume that $|S_1|<|S_2|$.
It is certain that $p=9$ or 11, and $|S_1|\leq 4$. Since $|S_1|=3$
provides $C[u,v]+uv$ is a 5-cycle, and $|S_1|=4$ provides
$C[v,u]+uv$ is a $(p-4)$-cycle, we assume that $|S_1|=1$ or 2.

If $|S_1|=1$, say $S_1=\{w\}$, then $uvwu$ bounds a 3-face by
Lemma \ref{clm0}. Let $G'$ be the graph obtained from $G-w$ by
inserting a new vertex into $uv$. Then, $G'\in {\cal G}$,
$\sigma(G')=\sigma(G)-1$, {\bf and the unbounded face of $G'$ is a
special face of $G'$ if $p=11$ since $f_o$ is one of $G$}. We can
extend $\phi$ to a 3-coloring $\phi'$ of $G'$. This produces a
contradiction because $\phi'$ and $\phi(w)$ yield a 3-coloring of
$G$ that extends $\phi$.

Assume $|S_1|=2$. Since $C[v,u]+uv$ is a $(p-2)$-cycle, and since
$G$ has neither adjacent triangles nor 5-cycles, $p=11$ and there
exists a 3-face sharing a unique edge with $f_o$ on $C[v,u]$. So,
$C[v,u]+uv$ is a separating 9-cycle, a contradiction to  Lemma
\ref{clm0}. $\hfill \rsq$

\medskip

\renewcommand{\baselinestretch}{1}
\begin{lemma}\label{11-cycles}
$N(u)\cap N(v)\cap int(C_1)=\emptyset$ for separating $11$-cycle
$C_1$ and $uv\in E(C_1)$.
\end{lemma}
\renewcommand{\baselinestretch}{2}
{\bf Proof.} Assume to the contrary that $x\in N(u)\cap N(v)\cap
int(C_1)$. By Lemma \ref{clm0}, $xuvx$ bounds a 3-face. We will
show that $C_1$ has neither claw-center nor d-claw-center. Then,
$C_1$ is a special cycle that contradicts Lemma \ref{clm0}.

{\bf Let $G'=G\setminus ext(C_1)$, and let $f'$ be the unbounded
face of $G'$. For each collapse $f''$ of $f'$, $xuvx$ is always
adjacent to $f''$, and a claw-center (resp. d-claw-center) of
$C_1$ is also one of $b(f'')$. We may assume that each claw-center
(resp. d-claw-center) of $C_1$ has three neighbors (resp. four
neighbors) on $C_1$}.

If $xw\in E(G)$ for some $w\in C_1\setminus\{u,v\}$, assume that
$u,v$ and $w$ clockwisely lie on $C_1$, then $|V(C_1(v,w))|\geq 5$
and $|V(C_1(w,u))|\geq 5$ since $G\in {\cal G}$, and hence
$|C_1|\geq 13$, a contradiction. If a vertex $y\in
int(C_1)\setminus\{x\}$ has three neighbors $z_1, z_2$ and $z_3$
on $C_1$, then by simply counting the number of vertices in
$C_1\setminus\{z_1, z_2, z_3\}$, $G$ must contain a 9-cycle $C_2$
with $x\in int(C_2)$, a contradiction to Lemma \ref{clm0} because
$C_2$ is a separating 9-cycle.

Assume that $\{a, b\}$ is a d-claw-center of $C_1$. Since $G$ has
no adjacent triangles, $|(N(a)\cup N(b))\cap C_1|\geq 3$. If
$(N(a)\cup N(b))\cap C_1$ has exactly three vertices, say $a_1,
a_2$ and $a_3$ clockwisely on $C_1$, we may assume that $a_1\in
N(a)\cap N(b)$, then $|V(C_1(a_1,a_2))|\geq 5$ and
$|V(C_1(a_3,a_1))|\geq 5$ that provide $|C_1|\geq 13$. So, assume
that $a$ has two neighbors $a_1, a_2\in C_1$, $b$ has two
neighbors $b_1, b_2\in C_1\setminus\{a_1, a_2\}$, and assume these
four vertices clockwisely lie on $C_1$.

If $a_1a_2\in E(C_1)$, then $|V(C_1(a_2, b_1))|\geq 4$ and
$|V(C_1(b_2,a_1))|\geq 4$ providing $|C_1|\geq 12$, a
contradiction. So, we may assume that $a_1a_2\not\in E(C_1)$ and
$b_1b_2\not\in E(C_1)$, i.e., $|V(C_1(a_1,a_2))|\geq 1$ and
$|V(C_1(b_1,b_2))|\geq 1$. By symmetry, we assume $x\in
int(C_1[a_1, b_1]\cup a_1abb_1)$.  By simply counting the number
of vertices in  $C_1\setminus\{a_1,a_2,b_1,b_2\}$, we get
$|C_1|>11$, a contradiction. $\hfill \rsq$

\medskip

\renewcommand{\baselinestretch}{1}
\begin{lemma}\label{length-2}
For $u,v\in C$ and $x\not\in C$, if $xu, xv\in E(G)$, then $uv\in
E(C)$.
\end{lemma}
\renewcommand{\baselinestretch}{2}
{\bf Proof.} Assume to the contrary that $uv\not\in E(C)$. By
Lemma~\ref{chordless}, $uv\not\in E(G)$. Let
$|V(C[u,v])|=l<|V(C[v,u])|$. Then, $3\leq l\leq {p+1\over 2}\leq
6$.

Since $C[u,v]\cup vxu$ is an $(l+1)$-cycle and $C[v,u]\cup uxv$ is
a $(p-l+3)$-cycle, $l\not\in \{4, 6\}$, and $l\neq 5$ whenever
$p=9$. If $l=5$ and $p=11$, $C[v,u]\cup uxv$ must bound a 9-face
by Lemma \ref{clm0}, then $f_o$ has to be adjacent to a 3-face
$f_1$ on $C[u,v]$, and hence $C[u,v]\cup vxu\cup b(f_1)$ yields a
7-cycle. So, $l=3$. Let $C[u,v]=uwv$.

If $p=11$, then there exists a 3-face sharing a unique edge with
$f_o$ on $C[v,u]$ that contradicts Lemma \ref{11-cycles} because
$C[v,u]\cup vxu$ is a separating 11-cycle. Therefore, $p=9$ and
$C[v,u]\cup vxu$ bounds a 9-face by Lemmas~\ref{clm0} and
\ref{chordless}. Let $G'$ be the graph obtained from $G\setminus
V(C(v,u))$ by inserting $5$ new vertices into $ux$. Then, $G'\in
{\cal G}$, $\sigma(G')<\sigma(G)$, and the unbounded face of $G'$
has degree $9$. We can extend $\phi(u), \phi(w)$ and $\phi(v)$ to
a 3-coloring $\phi'$ of $G'$ with $\phi'(u)\neq \phi'(x)$. But
$\phi'$ and $\phi$ yield a 3-coloring of $G$ that extends $\phi$,
a contradiction. $\hfill \rsq$

\medskip

\renewcommand{\baselinestretch}{1}
\begin{lemma}\label{4-face}
$G$ contains neither $4$-cycles nor $6$-cycles.
\end{lemma}\renewcommand{\baselinestretch}{2}
{\bf Proof.} First assume to the contrary that $G$ contains a
4-cycle. Assume that $C_1$ is a separating 4-cycle. Let $\psi$ be
an extension of $\phi$ on $G\setminus int(C_1)$, and let $G_1$ be
the graph obtained from $G\setminus ext(C_1)$ by inserting five
new vertices into an edge of $C_1$. If $p\neq 3$ then $|C\setminus
C_1|\geq 6$ since $C$ is chordless, and hence $|ext(C_1)|\geq 6$.
If $p=3$ then $|C\cap C_1|\leq 1$ and hence $E(C)\cap
E(C_1)=\emptyset$, again $|ext(C_1)|\geq 6$ because every face
incident with some edge on $C_1$ is a $4^+$-face. Therefore,
$\sigma(G_1)<\sigma(G)$, and we can extend the restriction of
$\psi$ on $C_1$ to $G_1$, and thus get a 3-coloring of $G$ that
extends $\phi$. So, we assume that $G$ contains no separating
4-cycles. We proceed to show that one can identify a pair of
diagonal vertices of a 4-cycle such that $\phi$ can be extended to
a 3-coloring of the resulting graph $G'$. Since any 3-coloring of
$G'$ offers a 3-coloring of $G$, this contradiction guarantees the
nonexistence of 4-cycles in $G$.

Let $f$ be an arbitrary 4-face of $G$ with $b(f)=uvwxu$.  If
$f\not\in N(f_o)$, $b(f)$ contains a pair of diagonal vertices
that are not on $C$. By symmetry, we assume that $u, w\in
b(f)\setminus C$ whenever $f\not\in N(f_o)$. Let $G_{u,w}$ be the
graph obtained from $G$ by identifying $u$ and $w$, and let
$r_{uw}$ be the new vertex obtained by identifying $u$ and $w$. It
is clear that $G_{u,w}$ contains no adjacent triangles since no
edge of $b(f)$ is contained in triangles. If $f\not\in N(f_o)$, it
is certain that $\phi$ is still a proper coloring of $C$ in
$G_{u,w}$. If $f\in N(f_o)$, we may assume that $u\in C$, then
$w\not\in C$ and $N(w)\cap C\subset \{x,v\}$ by Lemmas
\ref{chordless} and \ref{length-2}, and thus $\phi$ is also a
proper coloring of $C$ in $G_{u,w}$ by letting
$\phi(r_{u,w})=\phi(u)$.

Since a cycle of length 5 or 7 in $G_{u,w}$ yields a 7-cycle or a
separating 9-cycle in $G$, $G_{u,w}\in {\cal G}$. Now we need only
to check that $f_o$ is still a special face in $G_{u,w}$ in case
of $p=11$. Assume that $p=11$.

{\bf We first consider the case that $N(f_o)$ has 4-faces. Choose
$f$ to be a 4-face in $N(f_o)$. By symmetry, we assume that $ux\in
E(C)$. Let $x_1x_2uxx_3$ be a segment on $C$. Since $f_o$ is
adjacent to a 3-face and has no ears, we may suppose that
$v\not\in C$ and $xx_3$ is not on 4-cycles. Assume that $N(w)\cap
N(x_1)$ has a vertex, say $w'$. $w'\not\in C$ by
Lemmas~\ref{chordless} and \ref{length-2}, and so $(C\cup
x_1w'wx)\setminus\{u,x_2\}$ is an 11-cycle. Let $f'$ be a 3-face
sharing a unique edge with $f_o$. Either $b(f')\cap \{x_1x_2,
x_2u\}\neq \emptyset$ produces a 7-cycle, or $b(f')\cap
(C\setminus \{u,x_2\})\neq \emptyset$ contradicts
Lemma~\ref{11-cycles}. So, $N(w)\cap N(x_1)=\emptyset$, and $f_o$
has no ears in $G_{u,w}$}.

If $C$ has a claw-center $z$, then $z$ has three neighbors on $C$.
Let  $y_1, y_2$ and $y_3$ be three neighbors of $z$ clockwisely on
$C$ in $G_{u,w}$. Then $y_i=r_{uw}$ for an $i$. Assume
$y_1=r_{uw}$. It is clear that $x\not\in \{y_2, y_3\}$, and
$y_2y_3\in E(C)$ by Lemma \ref{length-2}. If $|V(C(x,y_2))|\leq
3$, then in $G$, $C(x,y_2)\cup xwzy_2\cup zy_3$ contains a cycle
of length 5 or 7. If $|V(C(y_3,u))|\leq 3$, then in $G$,
$C(y_3,u)\cup C_1\cup wzy_2\cup zy_3$ contains a cycle of length 5
or 7, or a separating 9-cycle. Therefore, $|V(C(x,y_2))|\geq 4$,
$|V(C(y_3,u))|\geq 4$, and hence $p\geq 12$, a contradiction.

Assume that $C$ has a d-claw-center $\{z_1,z_2\}$ in $G_{u,w}$.
Since $C$ has no claw-center in $G_{u,w}$, $|N(z_i)\cap C|=2$,
$i=1,2$. Let $N(z_1)\cap C=\{y_1,y_2\}$ and $N(z_2)\cap
C=\{y_3,y_4\}$. Since $G$ contains no adjacent triangles,
$\{y_1,y_2\}\cap \{y_3,y_4\}=\emptyset$ by Lemma \ref{length-2}.
Since $f_o$ is a special face in $G$, we may assume that
$y_2=r_{uw}$. Then, $y_3y_4\in E(C)$ by Lemma \ref{length-2}.
Using the similar argument as used in the last paragraph, we get
$p\geq 12$ by counting the number of vertices in $C(x,y_3),
C(y_4,y_1)$ and $C(y_1,u)$, a contradiction.

{\bf Suppose that $N(f_o)$ has no 4-faces. If $b(f)\cap
C\neq\emptyset$, both  $u$ and $w$ have no neighbor on
$C\setminus\{v,x\}$ by Lemmas~\ref{chordless} and \ref{length-2}.
If every 4-face shares no common vertex with $f_o$, we may suppose
that $w$ has no neighbor on $C$. In either case, it is
straightforward to check that $f_o$ has no ears in $G_{u,w}$}. $C$
has a claw-center $z$ provides $z=r_{u,w}$, and $C$ has a
d-claw-center provides $r_{u,w}$ is in the d-claw-center. In
either case, one may get a contradiction that $p\geq 12$ by almost
the same arguments as above.

\bigskip

Now, assume that $C'$ is a 6-cycle of $G$. Since $G$ contains no
4-cycles as just proved above, every face incident with some edge
on $C'$ is a $6^+$-face. If $C'$ is a separating cycle, it is not
difficult to verify that $|ext(C')|\geq 4$, then by letting $G''$
be the graph obtained from $G\setminus int(C')$ by inserting three
vertices into an edge of $C'$, we can first extend $\phi$ to
$G\setminus int(C')$, and then extend the restriction of $\phi$ on
$C'$ to $G''$, and thus get an extension of $\phi$ on $G$. So, we
assume that $G$ has no separating 6-cycles.


Let $f'$ be an arbitrary 6-face. If $b(f')\cap C\neq \emptyset$,
we choose $u_0$ to be a vertex in $b(f')\cap C$, and choose $u_1$
to be a vertex in $b(f')\setminus C$. If $b(f')\cap C=\emptyset$,
since $G$ contains no $l$-cycle for $l=4,5$ or 7, there must be a
vertex on $b(f')$ that has no neighbors on $C$, we choose such a
vertex as $u_1$. Let $b(f')=u_0u_1\ldots u_5u_0$, and let $H$ be
the graph obtained from $G$ by identifying $u_1$ and $u_5$, $u_2$
and $u_4$, respectively. Since $H$ contains no adjacent triangles,
and any 5-cycle (7-cycle) of $H$ yields a 7-cycle (separating
9-cycle) in $G$, $H\in {\cal G}$.

We will show that $\phi$ is still a coloring of $f_o$ in $H$. It
is trivial if $b(f')\cap C=\emptyset$, since the operation from
$G$ to $H$ is independent of $\phi$. Assume that $b(f')\cap C\neq
\emptyset$. Then, $u_0\in C$ and $u_1\not\in C$ by our choice, and
$u_2\not\in C$ and $N(u_1)\cap C=\{u_0\}$ by Lemma \ref{length-2}.
If either $u_2$ has no neighbors on $C$, or $u_4\not\in C$, then
we are done. Otherwise, assume that $u_4\in C$ and $u_2$ has a
neighbor, say $z$, on $C$, and assume that $u_0, z$ and $u_4$ lie
on $C$ clockwisely. Since $G$ contains no 5-cycles, $u_0u_4\not\in
E(G)$, and hence $u_5\in C$ by Lemma \ref{length-2}. Since $G$
contains no 4-cycles and no separating 6-cycles, $|V(C(u_0,
z))|\geq 4$, $|V(C(z,u_4))|\geq 4$, and hence $p\geq 12$, a
contradiction.

Finally, we will prove that $f_o$ is still a special face in $H$
in case of $p=11$. Then, a contradiction occurs again since $\phi$
can be extended to $H$ that offers an extension of $\phi$ to $G$,
this will end the proof of Lemma \ref{4-face} and also the proof
of our theorem.

The proof technique is again, as used repeatedly, to derive a
contradiction by counting the number of vertices on the segments
divided by the vertices adjacent to some claw-center or
d-claw-center of $C$. We leave the case that $b(f')\cap C\neq
\emptyset$ to the readres, and proceed only with the case
$b(f')\cap C=\emptyset$. {\bf Suppose that every 6-face has no
common vertex with $f_o$. Note that the above procedure holds for
an arbitrary 6-face of $G$, and note that $G$ has neither 4-cycles
nor separating 6-cycles as just proved, it is straightforward to
check that we can choose $f'$ to be a 6-face such that $f_o$ has
no ears in $H$}. Assume that $p=11$ but $f_o$ is not a special
face in $H$. Let $r_{1,5}$ and $r_{2,4}$ be the vertices obtained
by identifying $u_1$ and $u_5$, and $u_2$ and $u_4$, respectively.


Assume that $C$ has a claw-center $y$ with three neighbors
$y_1,y_2$ and $y_3$, clockwisely on $C$ in $H$. By symmetry, we
may assume that $y=r_{1,5}$, and assume that $y_1u_1\in E(G)$ and
$y_2u_5, y_3u_5\in E(G)$. Then, $y_2y_3\in E(C)$ by Lemma
\ref{length-2}. Since $G$ contains no adjacent triangles, contains
no cycles of length 4,5 and 7, and contains no separating
9-cycles, $|V(C(y_1,y_2))|\geq 4$, $|V(C(y_3,y_1))|\geq 5$, and
hence $p\geq 12$, a contradiction.

Assume that $C$ has a d-claw-center $\{z_1,z_2\}$ in $H$. Then,
each of $z_1$ and $z_2$ has two neighbors on $C$ and these four
vertices are all distinct. By symmetry, we may assume that
$z_1=r_{1,5}$. $z_2$ may be $u_0$, $r_{2,4}$ or a vertex not on
$C\cup C'$. In each case, the same argument as above ensures that
$p\geq 12$. This contradiction completes the proof of Lemma
\ref{4-face}. $\hfill \rsq$

\vskip 8pt

Our proof is then completed because by the assumption in Theorem
\ref{th-1}, $G$ contains either 4-cycles or 6-cycles. $\rsq$

\vskip 10pt

\renewcommand{\baselinestretch}{0.8}

\end{document}